\documentclass[10pt]{article}
\usepackage{amssymb,pst-node}
\def\sym{\mathbb}
\def\OV#1{\overline{#1}}

\def\B{{\sym B}}\def\M{{\sym M}}

\def\N{{\sym N}}

\def\Z{{\sym Z}}

\def\cqfd{\hspace{3.5mm} \hfill \vbox{\hrule height 3pt depth 2 pt width
2mm}
\bigskip}

\newtheorem{theorem}{Theorem}
\newtheorem{proposition}[theorem]{Proposition}
\newtheorem{lemma}[theorem]{Lemma}

\newtheorem{definition}[theorem]{Definition}
\newtheorem{note}{Note}
\newtheorem{remark}[note]{Remark}
\newtheorem{example}[note]{Example}

\newcommand{\shuffle}{
\kern 1pt
\rule{.3pt}{4pt}\rule{4.5pt}{.3pt}\rule{.3pt}{4pt}\rule{4.5pt}{.3pt}\rule{.3
pt}{4pt}
\kern 1pt
}

\begin{document}
\title{\bf Congruences Compatible with the Shuffle Product}
\author{G\'erard Duchamp and Jean-Gabriel Luque\\
LIFAR, Facult\'e des Sciences et des Techniques,\\
 76821Mont-Saint-Aignan CEDEX, France. }

\maketitle
\begin{abstract}
This article is devoted to the study of monoids which can be endowed with a shuffle product with coefficients in a 
semiring. We 
show that, when the multiplicities do not belong to a ring with  prime characteristic, such a 
monoid is  a monoid of traces. When the characteristic is prime, we give  a 
decomposition of the congruences $\equiv$ (or relators $R$) such that $A^*/_\equiv=\langle A; 
R\rangle$ admits a 
shuffle product. This decomposition involves only addition of primitive elements to the successive 
quotients. To end with, we study the compatibility with Magnus transformation and  
examine the case of  congruences which are homogeneous for some weight function. The existence of such a weight 
function is also showed for congruences of depth one.
\end{abstract}
\section{Introduction}

$\quad$Partially commutative structures share many nice combinatorial
properties
with
their free counterpart \cite{DK2,Schm,Vi3}. They have also been use
intensively in the computer science area \cite{Ga2,GaMa1}.
These structures can be rapidly described as being presented with
$$<A;\ \{[a,b]\}_{(a,b)\in \theta}>_{\bf cat}$$
where $A$ is a set of generators (an alphabet), $\theta\in A^2$ an
unoriented
graph without loop, {\bf cat} a suitable category (Monoid, Group,
$K$-associative algebras, Lie algebras) and $[a,b]$
expresses the fact that $a$ and $b$ commute in the structures of {\bf cat}.

The simplest (and already bearing all the combinatorial power \cite{DL}) of
these structures is the partially commutative monoid
$<A;\ \{ab=ba\}_{(a,b)\in \theta}>_{\bf Mon}$ (denoted $\M(A,\theta)$)
whose algebra over a semiring $K$ ($K<A,\theta>=K[\M(A,\theta)]$) is called
the algebra of partially commutative polynomials.\\
It is well known that, if a congruence $\equiv$ (i.e. an equivalence over
the free monoid $A^*$) is generated by commutations (i.e. is the kernel of a
natural morphism $nat: A^*\rightarrow \M(A,\theta)$) then the
shuffle product of two classes is a class-sum \footnote{The structure
constants are the partially commutative subword coefficients \cite{DR,Schm}.}
(This
can be verified by hand in the general case and in case $K$ is a ring, is
related to the fact that $K<A,\theta>$ is the envelopping algebra of the
- free - Lie algebra generated by the letters \cite{Th}).

This amounts to say that, if two series
$S_i=\sum_{w\in A^*}(S_i,w)w;\ i=1,2$ are $\equiv$-saturated (i.e. are
constant on every class of $\equiv$) then their shuffle product
$S_1\shuffle S_2$ is again saturated, which in turn reflects, by duality,
the
fact that there exists a coproduct $c_\equiv: K[A^*/_\equiv]\rightarrow
K[A^*/_\equiv]\otimes
K[A^*/_\equiv]$ such that the square
\begin{equation}\label{EQSH1}
\begin{array}{c@{\hskip 1cm}c}
\rnode{a}{K<A>}& \rnode{b}{K<A>\otimes K<A>} \\[1cm]
\rnode{c}{K[A^*{/_\equiv}]} &
\rnode{d}{K[A^*{/_\equiv}]\otimes K[A^*{/_\equiv}]}
\end{array}
\ncline{->}{a}{b}\Aput{c}
\ncline{->}{b}{d}\Aput{nat\otimes nat}
\ncline{->}{c}{d}\Aput{c_\equiv}
\ncline{->}{a}{c}\Bput{nat}
\end{equation}
is commutative ($c$ being the coproduct dual to the shuffle).
The converse was shown by {\sc Duchamp} and {\sc Krob} \cite{factor} with
the
restriction that $K$ be a ring of characteristic $0$.

The aim of our paper
is to discuss the property defined by (\ref{EQSH1}) which will be called,
throughout the paper $K-\shuffle$ compatibility.\\
Let us mention here that only the shuffle product is worth and gives
rise to such an interesting discussion as the compatibility with the other
classical rational laws (sum, external product, star and also Cauchy,
Hadamard and infiltration products \cite{lncs}) over series give weaker or
equivalent
results (for full details see Commentary \ref{Com}).
\smallskip
\noindent The structure of the paper is the following.

First (section 2), we prove that the result of {\sc Duchamp} and {\sc Krob}
holds in
almost every case, the only exception being the rings of prime
characteristic
for which differences of words that are primitive elements can occur. The
iteration of these adjunctions exhausts every finitely generated congruence.
More formally, we have:
\smallskip
\begin{theorem}\label{T1}
Let $\equiv$ be a finitely generated congruence on $A^*$ which is
$K-\shuffle$ compatible then
\begin{enumerate}
\item\label{IT1} If $K$ is not a ring or if the characteristic of $K$
($ch(K)$) is not prime then $A^*/_\equiv$ is a partially commutative monoid.
\item\label{IT2} If $ch(K)=p$ is prime, then a finite $R\subset A^*\times
A^*$ exists with a partition $R=\bigcup_{i\in[1,n]}R_i$ such that for $i\in[1,n]$, $R_i$ consists 
in
pairs
$(u,v)\in R$ such that (the image of) $u-v$ is primitive in
$K[A^*/_{\equiv_{\bigcup_{j=1}^{i-1}R_j}}].$
\end{enumerate}
\end{theorem}
\smallskip

In a second part (section 3), we investigate some properties of the new
family of
monoids. We prove that these monoids admit a Magnus transformation and we use it to examine the 
properties of cancellability, gradation, embeddability in a group and roots.
\begin{theorem}\label{T2}
Let $\equiv$ be a finitely generated congruence on $A^*$ which is
$K-\shuffle$ compatible and generated by relators of the form $u\equiv v$
where $u-v$ is primitive. We have the following alternatives.
\begin{enumerate}
\item The monoid $A^*/_\equiv$ is not cancellable.
\item A weight function exists $\omega:A\rightarrow {\sym N}^*$ for
which $\equiv$ is homogeneous and $A^*/_\equiv$ embeds in a group.
\end{enumerate}
\end{theorem}
\section{Compatibility with the Shuffle Product}
This section is devoted to study the $K-\shuffle$ compatibility of
a congruence. First, we address the problem in the general case ($K$ is a
semiring) and give some results common to all the semirings. We treat the
case when $K$ is the boolean semiring and when it is a ring, we develop some
results about the primitive polynomials which are useful to prove our
Theorem
\ref{T1}. In the last paragraph of this section, we sketch the proof.
\subsection{General Properties}
In the sequel we will denote $\N.1_K$ the subsemiring of a semiring $K$
generated by $1_K$.\\
We first remark that a congruence is $K-\shuffle$ compatible if and only
if it is $\N.1_K-\shuffle$ compatible. The property below is straightforward
from this observation.
\begin{lemma}\label{LETPHI}
Let $\phi:K_1\rightarrow K_2$ be a morphism of semirings then
\begin{enumerate}
\item If $\equiv$ is $K_1-\shuffle$ compatible then it is $K_2-\shuffle$
compatible.
\item If $\phi$ is into the converse holds.
\end{enumerate}
\end{lemma}
This shows that, in order to study the $K-\shuffle$ compatibility of a
congruence, it
suffices to study its $\N.1_K-\shuffle$ compatibility.\\
\begin{remark}
If $\equiv$ is $\N-\shuffle$ compatible then it is $K-\shuffle$ compatible
for each semiring $K$.
\end{remark}
The properties below will be useful in the sequel, their proof are easy and
we omit them.
\begin{lemma}\label{LProp}
\begin{enumerate}
\item If $\equiv_1$ and $\equiv_2$ are $K-\shuffle$ compatible
congruences then $\equiv_1\vee\equiv_2$ and $\equiv_1\wedge\equiv_2$ also are (supremum and
infimum are defined with respect to relation "is coarser than").
\item Let $R$ be a relator on $A^*$. The congruence $\equiv_R$ generated by
$R$ is $K-\shuffle$ compatible if and only if for each pair $(w_1,w_2)\in R$ we have
$c(w_1)\equiv_R^{\otimes 2}c(w_2)$.
\item Let $u\in A^+$ and $n$ be the maximal integer such that $u$ can be
written as $u=u_1a^n$ with $a\in A$. Then for each $K$ we have $(c(u),u_1\otimes a^n)=1$.
\item Each congruence generated by relators of the form $a\equiv b$ (LI)
or $dc\equiv cd$ (LC) with $a,b,c,d\in A$ is $K-\shuffle$ compatible.
\item Let $B\subseteq A$ be a subalphabet of $A$. If $\equiv$ is
$K-\shuffle$ compatible then its restriction to $B^*$ is.
\end{enumerate}
\end{lemma}
\subsection{The Boolean Case}
This paragraph deals with the case when $K=\B$ (the boolean semiring). The
$\B-\shuffle$ compatibility is completely characterized by the following
proposition.
\begin{proposition}\label{PBool}
A congruence is $\B-\shuffle$ compatible if and only if it is generated by
relators like $a\equiv 1$ (LE), $a\equiv b$ (LI) or $ab\equiv ba$ (LC), with $a,b\in A$.
\end{proposition}
{\em Sketch of the proof.} The "if" part is immediate using Properties (1),
(2) and (4) of Lemma \ref{LProp}.\\
For the converse, one defines a suitable section $S\subset A$ such that 
$S^*\displaystyle\mathop\rightarrow^{nat_S}A^*/_\equiv$ is onto and $\equiv_S=\mbox{Ker }(nat_S)$ 
is generated by only (LC) relators. We have successively verify that $\equiv_S$ is 
multihomogeneous and then, in fact, a partially commutative congruence.\cqfd
\subsection{Primitive Elements}
We suppose now that $K$ is a ring and we examine the link between
$K-\shuffle$ compatibility and primitivity of polynomials.
\begin{definition}
Let $K$ be a ring and $\equiv$ be a $K-\shuffle$ compatible congruence. A
polynomial $P\in K[A^*/_\equiv]$ is called primitive if and only if $c_\equiv(P)=P\otimes
1+1\otimes P$.
\end{definition}
As usual we have the following property.
\begin{proposition}\label{Lie}
The submodule $Prim(K[A^*/_\equiv])$ of primitive polynomials
endowed with the Lie bracket $[\quad,\quad]$ is a Lie algebra.
\end{proposition}
We suppose now that $\equiv=\equiv_\theta$ is a relation generated by
commutations.
Recall first that the free partially commutative monoid is
$$\M(A,\theta)=<A|\{ab=ba\}_{(a,b)\in \theta}>.$$
The monoid $\M(A,\theta)$ can be totally ordered by a relation $<_{std}$
in the following way :
\[t<_{std} t'\Leftrightarrow std(t)<_{lex} std(t')\]
where $std(t)$ denotes the maximal word for the lexicographical order in the
commutation class $t$.
Using this order Lalonde \cite{La,KL} has generalized the notion of
Lyndon word : the set of Lyndon traces is defined as the set of connected
and primitive\footnote{In the sense of traces (i.e. a trace is primitive if
it can not be written as the power of an other trace)} traces minimal in
their
conjugate classes and denoted $Ly(A,\theta)$. In his thesis Lalonde has
shown the following theorem.\\ \\
{\bf Theorem} {\it (Lalonde){\bf .} Let $\Lambda: Ly(A,\theta)\rightarrow
L_K(A,\theta)$ (the - free - Lie algebra generated by the letters in
$K<A,\theta>$) be the mapping defined by
\[
\left\{
\begin{array}{ll}
\Lambda(a)=a&\mbox{ if } a\in A\\
\Lambda(l)=[\Lambda(l_1),\Lambda(l_2)]&\mbox{ if } l=l_1l_2,\ l_1, l_2\in
Ly(A,\theta)\mbox{ and }|l_2|\mbox{ minimal}
\end{array}\right.
\]
then $(\Lambda(l))_{l\in Ly(A,\theta)}$ is a basis of $L_K(A,\theta)$ as a
$K$ module and for each Lyndon trace $l\in Ly(A,\theta)$ one has
\[\Lambda(l)=l+\sum_{t>_{std}l}\beta_tt\]
}
We have the following proposition.
\begin{proposition}\label{U-V}
Let $p$ be a prime integer and $u-v$ be a primitive polynomial of
$\Z/p\Z<A,\theta>$ such that $u,v\in\M(A,\theta)$. We can prove either $u=a^{p^\alpha}$ and
$v=b^{p^\beta}$ or $u=a^{p^\alpha}b^{p^\beta}$ and
$v=b^{p^\beta}a^{p^\alpha}$.
\end{proposition}
{\em Sketch of the proof.}
We denote $L_{\Z/p\Z}^{(p)}(A,\theta)$ the Lie algebra of the primitive
polynomials. We can prove successively that
\begin{enumerate}
\item[a)]\label{E1} The only monomials $u$ which are
primitive are of the form $u=a^{p^\alpha}$.
\item[b)]\label{E2} If $a_1\cdots a_n-b_1\cdots b_m$ is a
primitive
polynomial then we have $a_1\cdots a_n=a_n\cdots a_1$ and $b_1\cdots
b_m=b_m\cdots b_1$ or $n=m$ and $a_1\cdots a_n=b_n\cdots b_1$\footnote{In
the sequel, the mirror image of a trace $u$ will be denoted by $\OV u$.}.
\item[c)]\label{E3} Using Lalonde's Theorem, the set
$\{(\Lambda(l))^{p^e}\}$ generates $L_{\Z/p\Z}^{(p)}(A,\theta)$ as a
$K$-module and we write $u-v$ in the form
\[u-v=l^{p^e}+\sum_{t>l^{p^e}}\gamma_tt.\]
Without restriction we can consider that $u<_{std}v$ and then $u=l^{p^e}$.
According to (b), we have to consider two cases
\begin{enumerate}
\item[i)] If $u=\OV u$ and $v=\OV v$, we prove that $u$ and $v$ are
primitive and
by (a) we get the claim.
\item[ii)] Suppose that $u=\OV v$. Remarking that, if a power $l^e$ ($e>0$)
is not of the form $a^\alpha$ or $a^\alpha b^\beta$, one has
$l^{e}=awb^\alpha c^\beta$ with $a,b,c\in A$, $w\in\M(A,\theta)$,
$\alpha,\beta>0$ with $\alpha+\beta$ maximal, $c\neq b$ and $a\neq c$. In
this case $(c(l^e),awc^\beta\otimes b^\alpha)=1$ (by Property (3) of Lemma
\ref{LProp}) and
$(c(\OV{l^e}),awc^\beta\otimes b^\alpha)=0$, which proves that if $u\neq
a^\alpha b^\beta$ and $u\neq a^\alpha$, $u-v$ is not primitive. Solving the
equations $c(a^\alpha-b^{\beta})=(a^\alpha-b^{\beta})\otimes
1+(a^\alpha-b^{\beta})\otimes1$ and $c(a^\alpha
b^\beta-b^{\beta}a^{\alpha})=(a^\alpha b^\beta-b^{\beta}a^{\alpha})\otimes
1+1\otimes (a^\alpha b^\beta-b^{\beta}a^{\alpha})$, we obtain that $\alpha$ and $\beta$ are
necessarily two powers of $p$.
\end{enumerate}
\end{enumerate}\cqfd
\subsection{Structure of Compatible Congruences}
This paragraph is devoted to sketch the proof of Theorem \ref{T1}. We
recall it here.\\ \\
{\bf Theorem \ref{T1}.} {\it
Let $\equiv$ be a $K-\shuffle$ compatible and finitely generated congruence
on $A^*$. Then
\begin{enumerate}
\item If $K$ is not a ring or if the characteristic of $K$
($ch(K)$) is not prime then $A^*/_\equiv$ is a partially commutative monoid.
\item If $ch(K)=p$ is prime, then a finite $R\subset A^*\times
A^*$ exists with a partition $R=\bigcup_{i\in[1,n]}R_i$ such that for $i\in[1,n]$, $R_i$ consists 
in
pairs
$(u,v)$ such that (the image of) $u-v$ is primitive in
$K[A^*/_{\equiv_{\bigcup_{j=1}^{i-1}R_j}}].$
\end{enumerate}
}
{\em Sketch of the proof.}
We prove first (\ref{IT2}). Suppose that $A$ is finite (one can restrict
ourselves to the letters of the words of the relators), let $K$ be a ring
and $\equiv$ be a finitely
generated congruence on $A^*$ which is $K-\shuffle$ compatible. It is always
possible to derive a finite set of relators $R\subset A^*\times A^*$ closed
in the following sense
$$u\equiv
v\mbox{ and }\max\{u,v\}<\max\{\max\{u',v'\}|(u',v')\in
R\}\mbox{ $\Longrightarrow$ }(u,v)\in R.$$
We construct $(R_i)_{i\in[0,n]}$ in the following way.
\begin{enumerate}
\item We set $S_0=R_0=\emptyset$.
\item For each $i>0$, $R_i$ is the set of the pairs $(u,v)\in
R-\bigcup_{j\leq i-1}S_j$ such that $u-v$ is primitive in
$K[A^*/_{\equiv_{\bigcup_{j\leq i-1}R_j}}]$.
\item The relator $S_i$ is the set of the pairs $(u,v)\in R-\bigcup_{j\leq
i-1}R_j$ such
that $u\equiv_{R_i}v$.
\end{enumerate}
One can prove that this process ends, remarking that for each $i$ if
$\equiv_{\bigcup_{j\leq i-1}R_j}\neq \equiv_R$ we have $R_i\neq\emptyset$.
Now let us prove (\ref{IT1}). We may consider two cases.
\begin{enumerate}
\item The semiring $K$ is not a ring or a ring of characteristic 0. If
$1_K+1_K=1_K$, as one has $\B\hookrightarrow K$ Lemma \ref{LETPHI} and
propositon \ref{PBool} prove that it is generated by (LE), (LI) or (LC)
relators (see Lemma \ref{PBool}).
If $1_K+1_K\neq 1_K$, Lemma \ref{LETPHI} implies that $\equiv$ is generated
by (LE), (LI) or (LC) relators , examining all these cases, we find that
only
(LE) is impossible. In the two remaining cases, $A^*/_\equiv$ is a free partially
commutative monoid.
\item The semiring $K$ is ring of characteristic $n\neq 0$ not prime. We
consider two cases.
\begin{enumerate}
\item We have $n\neq p^\alpha$ with $p$ prime and $\alpha>1$. At
least two prime factors $p_1$ and $p_2$ of $n$ exist. By Lemma \ref{LETPHI}
$\equiv$ is $\Z/p_1\Z-\shuffle$ compatible and $\Z/p_2\Z-\shuffle$
compatible. Proposition \ref{U-V} and assertion (\ref{IT2}) imply that
$\equiv$ is generated only by (LI) or (LC) relators.
\item We have $n=p^m$, then by Lemma \ref{LETPHI} $\equiv$ is $\Z/p^2\Z-\shuffle$
compatible. Again by Lemma \ref{LETPHI} it is $\Z/p\Z$-compatible which
implies, using Proposition \ref{U-V}, that the only primitive polynomials
$u-v$ are of the form $a^{p^\alpha}-b^{p^\beta}$ or
$a^{p^\alpha}b^{p^\beta}-b^{p^\beta}a^{p^\alpha}$. Remarking that
$\left(p^\alpha\atop p^{\alpha-1}\right)\neq 0\ [p^2]$, we find that these
relators occur only when $\alpha=\beta=0$. Using assertion (\ref{IT2}),
again
Proposition \ref{U-V} gives the result.
\end{enumerate}
\end{enumerate}
\cqfd\\ \\
Such a family $(R_i)_{i\in[1,n]}$ will be called a {\tt primitive
partition} of $\equiv$. The minimal length of the primitive partitions will
be called the {\tt depth} of $\equiv$.
\begin{example}
\begin{enumerate}
\item All the congruences generated by relators of the form
$a^{p^\alpha}=b^{p^\beta}$ or
$a^{p^\alpha}b^{p^\beta}=b^{p^\beta}a^{p^\alpha}$
are $\Z/p\Z-\shuffle$ compatible with depth 1.
\item The congruence generated by
$$
\left\{\begin{array}{l}a^{2^{\alpha}}b^{2^\beta}a^{2^\alpha}b^{2^\beta}=b^{2
^\beta}a^{2^\alpha}b^{2^\beta}a^{2^{\alpha}}\\
a^{2^{\alpha+1}}b^{2^\beta}=b^{2^\beta}a^{2^{\alpha+1}}\\
b^{2^{\beta+1}}a^{2^\alpha}=a^{2^\alpha}b^{2^{\beta+1}}\end{array}\right.$$
with $\alpha,\beta\in\N-\{0\}$ are $\Z/2\Z-\shuffle$ compatible with
depth\footnote{We only know congruences with depth $2$.} $2$.
\end{enumerate}
\end{example}
\section{Group Properties}
\subsection{Compatibility with Magnus transformation}
In this paragraph we show that, when $\equiv$ is a $K-\shuffle$ compatible
congruence ($K$ semiring), one can define a Magnus transformation on
$K[A^*/_\equiv]$ . Recall that the Magnus transformation is the unique
endomorphism of $K<A>$ such that $\mu(a)=1+a$ for each letter $a$.
One has here
\begin{lemma}\label{Magnus}
Let $K$ be a semiring and $\equiv$ be a congruence $K-\shuffle$ compatible.
Then it exists an unique morphism $\mu_\equiv$ such that the square
\begin{equation}
\begin{array}{c@{\hskip 1cm}c}
\rnode{a}{K<A>}& \rnode{b}{K<A>} \\[1cm]
\rnode{c}{K[A^*{/_\equiv}]} &
\rnode{d}{K[A^*{/_\equiv}]}
\end{array}
\ncline{->}{a}{b}\Aput{\mu}
\ncline{->}{b}{d}\Aput{nat}
\ncline{->}{c}{d}\Aput{\mu_\equiv}
\ncline{->}{a}{c}\Bput{nat}
\end{equation}
is commutative.
\end{lemma}
{\bf Proof} It suffices to remark that $\mu=(Id\otimes ev)\circ c$
where $ev$ is the linear mapping from $K<A>$ on $K$ sending each
word $w\in A^*$ on $1$. This application is constant over $A^*$,
then it is compatible with the congruence $\equiv$ (i.e. it exists
an application $ev_\equiv$ sending each class of word to 1).
Furthermore, as $\equiv$ is $K-\shuffle$ compatible, $c_\equiv$
exists and, therefore our morphism is $\mu_\equiv=(Id\otimes
ev_\equiv)\circ c_\equiv$.\cqfd
\begin{remark}\label{Com}
Before closing the general case, let us mention that the problem of
compatibility can also be addressed for other rational laws. This laws are
"$.$" (concatenation or Cauchy product), $*$ (star operation, partially
defined),$\times$ (external product) and $+$ (union or sum) for the first
kind (which is the realm of Kleene-Sch\"utzenberger Theorem \cite{Kle}) and
for the second kind $\odot$ (Hadamard product), $\shuffle$ (shuffle
product), $\uparrow$ (infiltration product \cite{Ochs}) and $\uparrow_q$
($q$-infiltration product\footnote{ It can be shown that $q$-infiltration
is the only dual law satisfying alphabetical and algebraic constraints
\cite{lncs}}, the dual of the coproduct $c_q=(a\otimes1+1\otimes
a)+qa\otimes a$). The results can be summarized as follows
\[
\begin{array}{|c|c|c|c|}
\hline
Kind&{\rm Laws}&\mbox{Compatible with commutation}&\mbox{Other}\\
\hline
&\times&\mbox{Yes}&\mbox{All}\\
&+&\mbox{Yes}&\mbox{All}\\
\mbox{First}&.&\mbox{No}&\mbox{?}\\
&$*$&\mbox{No}&\mbox{?}\\
\hline
&\odot&\mbox{Yes}&\mbox{All}\\
&\shuffle&\mbox{Yes}&\mbox{depends of $K$}\\
\mbox{Second}&\uparrow&\mbox{Yes}&\mbox{as }\shuffle\\
&\uparrow_q&\mbox{Yes}&\mbox{as }\shuffle\\
\hline
\end{array}
\]
\end{remark}
\subsection{Homogeneous Congruences Compatible with the Shuffle}
In this paragraph, we examine the congruences which are homogeneous for some
weight function
$\omega: A\rightarrow \N^+$. The function $\omega$ is extended to $A^*$ by
$\omega(w)=\sum_{a\in{\rm Alph}(w)}|w|_a\omega(a)$, where $|w|_a$ is the
partial degree of $w$ with respect to the letter $a$. Remarking that
$\omega$ is a
morphism the
following result is straigthforward.
\begin{lemma}\label{LETR}
Let $R\subset A^*\times A^*$. The following assertions are equivalent.
\begin{enumerate}
\item For each $(u,v)\in R$, $\omega(u)=\omega(v)$.
\item For each $u,v\in A^*$, $u\equiv_R v\Rightarrow \omega(u)=\omega(v)$.
\end{enumerate}
\end{lemma}
\begin{definition}
We will say that $\equiv$ is homogeneous (for $\omega$) if and only if it
satisfies the assertions of Lemma \ref{LETR}.
\end{definition}
{\bf Theorem \ref{T2}} {\it
Let $\equiv$ be a finitely generated congruence on $A^*$ which is
$K-\shuffle$ compatible and generated by relators of the form $u\equiv v$
where $u-v$ is primitive. We have the following altenative.
\begin{enumerate}
\item The monoid $A^*/_\equiv$ is not cancellable.
\item A weight function exists $\omega:A\rightarrow {\sym N}^*$ for
which $\equiv$ is homogeneous and $A^*/_\equiv$ embeds in a group.
\end{enumerate}
}
{\em Sketch of the proof.} Suppose that $A^*/_\equiv$ is cancellable. If
$ch(K)$ is not prime then $A^*/_\equiv$ is a free partially commutative
monoid, and the result is straightforward from \cite{DK3}. If $ch(K)$ is prime,
$\equiv$ is generated by relators of type
$a^{p^\alpha}b^{p^\beta}=b^{p^\beta}a^{p^\alpha}$ (pLC) and
$a^{p^\alpha}=b^{p^\beta}$ (pLI). As all the (pLC) relators are
multihomogeneous
it suffices to prove the existence of a weight function for which (pLI) is
homogeneous. We prove that, if we have two relators $a^{p^{\alpha_i}}\equiv
b^{p^{\beta_i}}$ ($i= 1,2$) one has
$\alpha_1-\beta_1=\alpha_2-\beta_2$. Denoting by $d_{a,b}$ this difference, the
result is a consequence of the following lemma whose proof is easy and left to
the reader. \begin{lemma}\label{LETG}
Let $G$ be the graph of an equivalence relation on $A$. Let
$d:G\rightarrow \Z$ be a function such that $$(a,b),(a,c)\in G\Rightarrow
d(a,b)+d(b,c)+d(c,a)=0.$$ Then,
\begin{enumerate}
\item It exists a (potential) function $h: A\rightarrow\Z$ such that
$d(a,b)=h(b)-h(a)$.
\item If $A$ is finite, we can choose $h$ positive.
\end{enumerate}
\end{lemma}
{\em End of the proof.}
We remark that $d_{a,b}+d_{b,c}+d_{c,a}=0$, according to Lemma \ref{LETG} we
can construct a function $h:A\rightarrow \N^+$ such that
$d_{a,b}=h(b)-h(a)$. With the weight function $\omega=p^h$, the congruence
$\equiv$ is homogeneous.
On the other hand we have $\mu_\equiv(A^*)\subset 1+{\cal M}_K(A^*/_\equiv)$
where
${\cal M}_K(A^*/_\equiv)$ is the ideal of series $S$ such that $(S,1)=0$.
But,
as $\equiv$ is homogeneous, $1+{\cal M}_K(A^*/_\equiv)$ is a group.
Furthermore
$\mu_\equiv(w)=w+\sum_{\omega(u)<\omega(w)}n_uu$, which implies that
$\mu_\equiv:
A^*/_\equiv\hookrightarrow 1+{\cal M}_K(A^*/_\equiv)$ is into. This proves the
result.\cqfd

We now give quickly some properties of the Magnus group
$1+{\cal M}_K(A^*/_\equiv)$.
It has been shown in \cite{DK3} that, if $ch(K)=0$, the function
$S\rightarrow S^n;\ n>0$ is one to one within the Magnus group (the monoid
then is partially commutative). We cannot expect such a property if $ch(K)=p$ because the 
congruence $a^p\equiv b^p$ can occur with
$a\not\equiv b$. However, we have
\begin{proposition} Let $K$ be a ring of prime characteristic $p$ and $\equiv$ a
$K-\shuffle$ compatible and cancellable congruence. Then:\\
i) If $q\not\equiv 0\ [p]$ the function $S\rightarrow S^q$ is one to one
within
the Magnus group $1+{\cal M}_K(A^*/_\equiv)$.\\
ii) If $q_i\not\equiv 0\ [p];\ i=1,2$ then $S^{q_1}T^{q_2}=T^{q_2}S^{q_1}$
implies $ST=TS$
\end{proposition}
{\bf Proof} Solving degree by degree
\begin{equation}
(1+\sum_{i=1}^\infty X_i)^q=(1+\sum_{i=1}^\infty Y_i)
\end{equation}
(with two auxilliary alphabets such that $deg(X_i)=deg(Y_i)$), it can easily
be seen that the series $\sum_{k=0}^\infty \left({\frac1q\atop k}\right)
X^k$ has its
coefficients in $\Z[1/q]$. This proves (i).\\
Now (ii) is a consequence of (i) as the hypothesis can be reformulated
$T^{-q_2}S^{q_1}T^{q_2}=S^{q_1}$.\\
\hfill\cqfd
\section{Conclusion}
We have seen that in ``almost every case'', the congruences compatible with
the shuffle product give partially commutative quotients. The only
degenerate
case occurs when the characteristic of the ground ring is prime and gives rise to a bunch of new 
phenomena. A process (primitive partitions) to
analyse these new congruences has been described. Moreover, for congruences
of depth
one, we have a complete description of the relators. This allows to prove,
thanks to a ``Magnus-type'' transformation, that the quotients are either not
cancellable or embeddable in a group. This transformation gives us the
opportunity to use some analytic tools as the series of the $q-root$ for $q$
prime to the characteristic. An infinite family of congruences of depth two
has been provided. The problem of describing higher depth remains however
open.


\begin{thebibliography}{ABC}
\bibitem{BR}
J. Berstel and C. Reutenauer,
{\em Rational Series and Their Languages,\/}
(EATCS Monographs on Theoretical Computer Science, Springer-Verlag Berlin,
1988).
\bibitem{Bo2}
N.Bourbaki,
{\em \'El\'ements de math\'ematiques, Groupes et alg\`ebres de Lie, Chap. 2
et 3\/}
(Hermann, Paris, 1972).
\bibitem{CF}
P.Cartier, D.Foata,
{\em Probl\`emes combinatoires de commutation et r\'earrangement\/},
Lect. Not. In Math., n 85, 1969.
%
\bibitem{CFL}
K.T.Chen, R.H.Fox, R.C.Lyndon,
{\em Free differential calculus IV- The quotient groups of the lower
central series,\/}
Ann. Of Math, 1958.
\bibitem{DR}
V. Diekert and G. Rozenberg,
{\em The book of traces\/}
(World Scientific, Singapour, 1995).
%
\bibitem{lncs}
G. Duchamp, M. Flouret, \'E. Laugerotte,
{\it Operations over Automata with Multiplicities,}
in Automato implementation procedeeding WIA, J.M. Champarnaud, D.Maurel and
D. Ziadi eds,{\bf 1660}, 183-191,1999.
%
\bibitem{factor} G. Duchamp , D. Krob, {\it Factorisations dans le mono\"\i
de
partiellement commutatif libre}, C.R. Acad. Sci. Paris, t. {\bf312},
s\'erie I (1991), 189-192.
%

\bibitem{DK2}
G.Duchamp, D.Krob,
{\em Free partially commutative structures\/}
{\em J. Algebra\/} {\bf 156}{-2} (1993) 318--361.
%
\bibitem{DK3}
G.Duchamp, D.Krob,
{\em Partially commutative Magnus transformation\/}
Int. J. of Alg. And comp.,
3-1, 1993,15-41.
%
%
\bibitem{DL}
G.Duchamp, J.G.Luque,
{\em Transitive Factorizations\/},
Colloque FPSAC'99 Barcelone,1999.
%
\bibitem{Fl}
M.Flouret,
{\em Contribution \`a l'algorithmique non commutative},
Th\`ese de doctorat, Univerit\'e de Rouen, 1999.
%
\bibitem{Ga2}
P. Gastin,
{\em Decidability of the Star problem in $A^*\times\{b\}^*$},
Information Processing Letters,
44,65-71,1992.
%
\bibitem{GaMa1}
S. Gaubert, J.Mairesse,
{\em Medeling and analysis of timed Petri nets using heap of pieces},
IEEE,Trans. Autom. Control, Vol 44, n4,683-697,1999.
%
\bibitem{Kle}
S.C. Kleene, {\it Representation of events in nerve nets and finite
automata}, Automata Studies, Princeton Univ. Press (1956), 3-42.
%
\bibitem{KL}
D. Krob and P.Lalonde,
{\em Partially commutative Lyndon words\/}
{\em Lect. Notes in Comput. Sci.\/} {\bf 665} (1993) 237--246.
%
\bibitem{La}
P.Lalonde,
{\em Contribution \`a l'\'etude des empilements}
(Th\`ese de doctorat, LACIM, 1991).
%
%
\bibitem{Lo}
M.Lothaire,
{\em Combinatorics on words\/},
Addison Wesley, 1983.
%
\bibitem{Ochs}
P. Ochsenschl\"ager, Binomialkoeffitzenten und Shuffle-Zahlen, Technischer
Bericht, Fachbereich Informatik, T. H. Darmstadt,1981.
\bibitem{Schm}
W.Schmitt,
{\em Hopf algebras and identities in free partially commutative monoids\/},
T.C.S. North Holland, 1990.
\bibitem{Th} J.Y. Thibon, {\it Int\'egrit\'e des alg\`ebres de s\'eries
formelles sur un alphabet
partiellement commutatif},T.C.S (1985), North-Holland.
\bibitem{Vi3}
X.G.Viennot,
{\em Heaps of pieces I: Basic definitions and combinatorial lemmas}
In G. Labelle et al., editors, Proceeding Combinatoire \'enumeratice,
Montr\'eal Quebec 1985,
nymber 1234 in Lectures notes in
Mathematics, 321-350, Berlin-Heidelberg-New York, 1986, Springer.
\end{thebibliography}
\end{document}